\documentclass[11pt,bezier]{article}
\pdfoutput=1
\usepackage{amsmath}
\usepackage{amsfonts,amsthm,amssymb}
\usepackage{exscale}
\usepackage{relsize}
\usepackage{graphicx}
\usepackage{epstopdf}
\usepackage{caption}
\usepackage{float}

\usepackage{amssymb}
\usepackage{amsmath}
\newtheorem{lemma}{Lemma}[section]
\newtheorem{theorem}[lemma]{Theorem}
\newtheorem{corollary}[lemma]{Corollary}

\usepackage{graphicx}
\makeatletter
\newcommand{\Rmnum}[1]{\expandafter\@slowromancap\romannumeral #1@}
\makeatother
\theoremstyle{definition}

\usepackage[justification=centering]{caption}
\begin{document}
\title{Extra Connectivity of Strong Product of Graphs\footnote{The research is supported by National Natural Science Foundation of China (11861066).}}
\author{Qinze Zhu, Yingzhi Tian\footnote{Corresponding author. E-mail: zqz@stu.xju.edu.cn (Q. Zhu); tianyzhxj@163.com (Y. Tian).} \\
{\small College of Mathematics and System Sciences, Xinjiang
University, Urumqi, Xinjiang 830046, PR China}}

\date{}
\maketitle

\begin{sloppypar}

\noindent{\bf Abstract } The $g$-$extra$ $connectivity$ $\kappa_{g}(G)$ of a connected graph $G$ is the minimum cardinality of a set of vertices, if it exists, whose deletion makes $G$ disconnected and leaves each remaining component with more than $g$ vertices, where $g$ is a non-negative integer. The $strong$ $product$ $G_1 \boxtimes G_2$ of graphs $G_1$ and $G_2$ is the graph with vertex set $V(G_1 \boxtimes G_2)=V(G_1)\times V(G_2)$, where two distinct vertices $(x_{1}, y_{1}),(x_{2}, y_{2}) \in V(G_1)\times V(G_2)$ are adjacent in $G_1 \boxtimes G_2$ if and only if $x_{1}=x_{2}$ and $y_{1} y_{2} \in E(G_2)$ or $y_{1}=y_{2}$ and $x_{1} x_{2} \in E(G_1)$ or $x_{1} x_{2} \in E(G_1)$ and $y_{1} y_{2} \in E(G_2)$. In this paper, we give the $g\ (\leq 3)$-$extra$ $connectivity$ of $G_1\boxtimes G_2$, where $G_i$ is a maximally connected $k_i\ (\geq 2)$-regular graph for $i=1,2$. As a byproduct, we get $g\ (\leq 3)$-$extra$ conditional fault-diagnosability of $G_1\boxtimes G_2$ under $PMC$ model.


\section{Introduction}
Throughout this paper, all the graphs are simple, that is, with neither loops nor multiple edges. We follow \cite{Bondy} for graph-theoretical terminology and notation which are not defined here. Let $G$ be a graph with vertex set $V(G)$ and edge set $E(G)$. For each vertex $u\in V(G)$, the neighbors $N_G(u)$ is the set of all vertices adjacent to $u$, and $d(u)$ is the degree of $u$. The closed neighborhood $N_G[u]$ of $u$ is $N_G(u)\cup\{u\}$. We denote by $\delta(G)$ the $minimum$ $degree$ of $G$. The graph $G$ is $k$-regular if $d(u)=k$ for any $u \in V(G)$. A vertex-cut in $G$ is a set of vertices whose deletion makes $G$ disconnected. It is a known fact that only complete graphs do not have vertex-cuts. The $connectivity$ $\kappa(G)$ of the graph $G$ is the minimum order of a vertex-cut in $G$ if $G$ is not a complete graph, and $\kappa(G)=|V(G)|-1$ otherwise. The girth of $G$, denoted by $girth(G)$, is the length of its shortest cycle in $G$ if it has any, otherwise, define $girth(G) = +\infty$.

Fault-tolerance is one of the main factors which should be considered in the design of an interconnection network. Usually, the topology structure of an interconnection network can be modeled by a graph $G$, where $V(G)$ represents the set of processors and $E(G)$ represents the set of links connecting processors in the network, the connectivity of $G$ is an important measurement for reliability and fault-tolerance of the network. In general, the larger $\kappa(G)$ is, the more reliable the network is. Because $\kappa(G) \leq \delta(G)$, a graph $G$ is said to be maximally connected if $\kappa(G)=\delta(G)$. Clearly, if $G$ is maximally connected and $d_{G}(u)=\delta(G)$, then $N_{G}(u)$ is a minimum vertex-cut and isolates the vertex $u$. Indeed, it is often expected that the system can work even if several of its elements are destroyed. Thus, a reliable design requires that the graph has maximum
connectivity. Besides, it is believed that the graph with the smallest number of minimum vertex-cuts is more reliable than others (see \cite{Boesch 1, Boesch 2}).

The connectivity of a graph $G$ has an obvious deficiency, that is, it tacitly assumes that all vertices adjacent to the same vertex of $G$ can potentially fail at the same time, which is almost impossible in the practical applications of networks. To overcome this deficiency, Harary \cite{Harary} proposed conditional connectivity. For a graph-theoretic property $\mathcal{P}$, the $conditional$ $connectivity$ $\kappa(G; \mathcal{P})$ is the minimum cardinality of a set of vertices whose deletion disconnects $G$ and every remaining component has property $\mathcal{P}$. Later, F{\`{a}}brega and Fiol \cite{Fabrega} introduced specific conditional connectivity named $g$-$extra$ $connectivity$. Let $g$ be a non-negative integer, a $g$-$extra$ $cut$ is a set of vertices in $G$ if it exists, whose deletion makes $G$ disconnected and leaves each remaining component with more than $g$ vertices. The $g$-$extra$ $connectivity$, denoted by $\kappa_{g}(G)$, is the minimum order of a $g$-$extra$ $cut$ if it has any, otherwise define $\kappa_g(G)= +\infty$. When we call a $g$-$extra$ $cut$ as a $\kappa_g$-cut in $G$ if and only if the order the $g$-$extra$ $cut$ is $\kappa_{g}(G)$. Under the definition of $g$-$extra$ $connectivity$, we define $\kappa(G)=\kappa_0(G)$ for any connected graph $G$ that is not a complete graph, so the $g$-$extra$ $connectivity$ can be seen as a generalization of the connectivity and more refined index for large systems and improves the accuracy of the reliability and fault tolerance of networks. The authors in \cite{Chang} pointed out that there is no polynomial-time algorithm for computing $\kappa_g$ for a general graph. Consequently, much of the work has been focused on the computing of the $g$-$extra$ $connectivity$ of some given graphs, see [3, 5, 8-10, 12-13, 15-17, 22-27, 29, 30] for examples.

The construction of new graphs from two given ones is routine. The method consists of joining together several copies of one graph according to the structure of another one, the latter being usually called the main graph of the construction. Since 1960 some relevant graph theory researchers have defined different types of graph products. The main difference between them comes from the number of inter-copy edges and the connection criterion. The most studied four standard graph products are the Cartesian product, the direct product, the strong product, and the lexicographic product.

The Cartesian product of two graphs $G_1$ and $G_2$, denoted by $G_1 \square G_2$, is defined on the vertex sets $V(G_1) \times V(G_2)$, and $(x_1, y_1)(x_2, y_2)$ is an edge in $G_1 \square G_2$ if and only if one of the following is true: ($i$) $x_1=x_2$ and $y_1 y_2 \in E(G_2)$; ($ii$) $y_1=y_2$ and $x_1 x_2 \in E(G_1)$. During the past years, many works focused on the $g$-$extra$ $connectivity$ of the Cartesian product of regular graphs. The connectivity of the Cartesian product of graphs was determined in \cite{Spacapan1}. Tian and Meng \cite{Tian} proved that $\kappa_{1}(K_{m} \square G_{2})=2 k_{2}+m-2$ for $m \geq 2$ and $\kappa_{1}(G_{1} \square G_{2})=2 k_{1}+2 k_{2}-2$, where $G_{i}$ is a $k_{i}(\geq 2)$-regular and maximally connected graph with $g(G_{i}) \geq 4$ for $i=1,2$. And Chen, Meng, Tian, and Liu \cite{Chen} gave that $\kappa_{2}(K_{m} \square G_{2})=3 k_{2}+m-3$ for $m \geq 3$ and $\kappa_{2}(G_{1} \square G_{2})=3 k_{1}+3 k_{2}-5$, where $G_{i}$ is a $k_{i}(\geq 2)$-regular and maximally connected graph with $g(G_{i}) \geq 5$ for $i=1,2$.

The strong product $G_1 \boxtimes G_2$ of $G_1$ and $G_2$ is the graph with the vertex set $V(G_1 \boxtimes G_2)=V(G_1) \times V(G_2)$, where two vertices $(x_{1}, y_{1})$, $(x_{2}, y_{2}) \in V(G_1) \times V(G_2)$ are adjacent in $G_1 \boxtimes G_2$ if and only if one of the following is true: ($i$) $x_1=x_2$ and $y_1 y_2 \in E(G_2)$; ($ii$) $y_1=y_2$ and $x_1 x_2 \in E(G_1)$; ($iii$) $x_1 x_2 \in E(G_1)$ and $y_1 y_2 \in E(G_2)$. Some properties regarding the minimum degree of $G_1\boxtimes G_2$ can be found in \cite{Sabidussi}. The structure of every miminimum vertex-cut in $G_1\boxtimes G_2$ was determined by {\v{S}}pacapan\cite{Spacapan}, as a corollary, he gave that $\kappa(G_1\boxtimes G_2)\geq$ min$\{\kappa(G_1)(1+\delta(G_2)),\kappa(G_2)(1+\delta(G_1))\}$. Zhu and Tian \cite{Zhu1} proved $\kappa_g(C_m\boxtimes C_n)=min\{2m, 2n, 2\lceil\, 2\sqrt{g+1}\,\rceil+4\}$ if $g\leq max\{n\lfloor\frac{m-2}{2}\rfloor-1, m\lfloor\frac{n-1}{2}\rfloor-1 \}$. And Casablanca, Cera, V{\'a}zquez, and Valenzuela \cite{Casablanca} gave a sufficient condition for the strong product of two maximally connected graphs, $G_1$ and $G_2$, to be maximally connected. These conditions are addressed in terms of the minimum degree and the girth of both $G_1$ and $G_2$.

\begin{theorem}\cite{Casablanca}\label{5}
Let $G_1$ and $G_2$ be two connected graphs with at least $3$ vertices and girth at least $4$. Then $G_1 \boxtimes G_2$ is maximally connected if both $G_1$ and $G_2$ are maximally connected and one of the following assertions holds:

(i) One graph has minimum degree $1$ and the other graph has girth at least $5$.

(ii) $\delta(G_i) \geq 2,\;i=1,2$.
\end{theorem}

Motivated by the results above, we will study the $g$-$extra$ $connectivity$ of the strong product of maximally connected graphs.

In the next section, we introduce some definitions and lemmas. In Section 3, we mainly argue the $g(\leq 3)$-$extra$ $connectivity$ of $G_1\boxtimes G_2$, where $G_i$ is a maximally connected $k_i\,(\geq 2)$-regular graph for $i=1,2$. The notation and result of $g\ (\leq 3)$-$extra$ conditional fault-diagnosability of $G_1\boxtimes G_2$ under $PMC$ model are given in Section 4. And the conclusion will be given in Section 5.

\section{Preliminary}
Let $V_1,\;V_2\subseteq V(G)$. Define $N_{V_1}(u)=\{\,v\in V_1|\, v \text{ is adjacent to u}\,\}$ is the set of neighbors of $u$ in $V_1$ and $N_{V_1}(V_2)=\{\,v\in V_1\backslash V_2|\, v \text{ is adjacent to a vertex }u\in V_2\}$ is the set of neighbors of $V_2$ in $V_1$. We denote $N_{V_1}[u]=N_{V_1}(u)\cup\{u\}$ and $N_{V_1}[V_2]=N_{V_1}(V_2)\cup V_2$. Sometimes we use a graph itself to represent its vertex set. For example, we use $N_{G_1}(V_2)$ instead of $N_{V_1}(V_2)$, where $G_1=G[V_1]$ is the induced subgraph of $G$ on $V_1$.

In the following, we present the same notations as \cite{Spacapan} related to the strong product graph. Let $G_1$ and $G_2$ be two graphs. Define two natural projections $p_1$ and $p_2$ on $V(G_1)\times V(G_2)$ as follows: $p_1(x,y)=x$ and $p_2(x,y)=y$ for any $(x,y)\in V(G_1)\times V(G_2)$. The subgraph induced by $\{(u, y)|u\in V(G_1)\}$ in $G_1\boxtimes G_2$, denoted by $G_{1y}$, is called a $G_1$-layer in $G_1\boxtimes G_2$ for each vertex $y\in V(G_2)$. Analogously, the subgraph induced by $\{(x, v)|v\in V(G_2)\}$ in $G_1\boxtimes G_2$, denoted by ${}_{x}G_2$, is called a $G_2$-layer in $G_1\boxtimes G_2$ for each vertex $x\in V(G_1)$. Clearly, a $G_1$-layer in $G_1\boxtimes G_2$ is isomorphic to $G_1$, and a $G_2$-layer in $G_1\boxtimes G_2$ is isomorphic to $G_2$.

Let $S\subseteq V(G_1\boxtimes G_2)$. For any $x\in V(G_1)$, denote $S\cap V({}_{x}G_2)$ by ${}_{x}S$, and analogously, for any $y\in V(G_2)$, denote $S\cap V(G_{1y})$ by $S_{y}$.

\begin{lemma}\cite{Zhu1}\label {0}
Let $g$ be a non-negative integer, and $G_i$ be a connected graph for $i=1,2$. Assume $G_1\boxtimes G_2$ has $g$-$extra$ cuts and $S$ is a $\kappa_g$-cut of $G_1\boxtimes G_2$.

(i) If ${}_{x}S\neq\emptyset$ for some $x\in V(G_1)$, then $|{}_{x}S|\geq \kappa(G_2)$.

(ii) If $S_{y}\neq \emptyset$ for some $y\in V(G_2)$, then $|S_{y}|\geq \kappa(G_1)$.

\end{lemma}

To show our main results, we give a lower bound of the number of vertices in a regular graph $G$ whose girth is at least 5.

\begin{lemma}\label{4}
If $G$ is a $k\ (\geq 2)$-regular graph with $girth(G)\geq 5$, then $|V(G)|\geq (girth(G)-2)(k-1)+2$ and the equality holds if and only if $k=2$.
\end{lemma}
\noindent{\bf Proof.} Since $d(u)\geq 2$ for all $u\in V(G)$, $G$ contains a cycle. Let $C$ be a cycle in $G$ with length of $girth(G)$. Denote $V(C)=\{u_1,u_2,\dots,u_{girth(G)}\}$. Then $(N_{G}(u_i)\backslash V(C))\cap (N_{G}(u_j)\backslash V(C))=\emptyset$ for any $i\neq j\in\{1,2,\dots,girth(G)\}$, otherwise, there is an another cycle whose length is less than $girth(G)$. Therefore $|V(G)|\geq |N_G[V(C)]|=\sum^{girth(G)}_{i=1}|N_{G}(u_i)\backslash V(C)|+|V(C)|=girth(G)(k-2)+girth(G)=girth(G)(k-1)\geq (girth(G)-2)(k-1)+2$, and the last equality holds if and only if $k=2$. The proof is thus complete. $\Box$

At the end of this section, we recall the famous Menger's Theorem \cite{Menger}.
\begin{theorem}(Menger's Theorem)\cite{Menger}\label{1}
Let $G$ be a connected graph with $x,y\in V(G)$. Then the minimum number of vertices separating vertex $x$ from vertex $y$ in $G$ is equal to the maximum number of internally disjoint $(x,y)$-paths in $G$ if $xy\notin E(G)$.
\end{theorem}

\section{The Extra Connectivity of the Strong Product of Two Regular and Maximally Connected Graphs}

Theorem $\ref{5}$ gives the sufficient condition for the strong product of two maximally connected graphs, $G_1$ and $G_2$, to be maximally connected. This implies that if $G_i$ is a $k_i\ (\geq 2)$-regular and maximally connected graph with $girth(G_i)\geq 5$ for $i=1,2$, then $\kappa(G_1\boxtimes G_2)=\delta(G_1\boxtimes G_2)=k_1k_2+k_1+k_2$.

\begin{theorem}\label{6}
Let $G_i$ be a $k_i\ (\geq 2)$-regular and maximally connected graph with $girth(G_i)\geq 5$ for $i=1,2$. Then $\kappa_{1}(G_1\boxtimes G_2)=$ min$\{2k_1k_2+2k_1-2,2k_1k_2+2k_2-2\}$.
\end{theorem}
\noindent{\bf Proof.} Denote $G=G_1\boxtimes G_2$. Let $u_1u_2$ be an edge in $G$, where $u_1=(x_1,y_1)$ and $u_2=(x_2,y_1)$. It is routine to verify that $N_{G}(\{u_1,u_2\})$ is a $1$-$extra$ $cut$ in $G$. Clearly, $|V(G-N_{G}[\{u_1,u_2\}])|> 2$, therefore if we can prove that $G-N_{G}[\{u_1,u_2\}]$ is connected, then $N_{G}(\{u_1,u_2\})$ is a $1$-$extra$ $cut$ in $G$. Denote $G'_1=G_1-N_{G_1}[\{x_1,x_2\}]$ and $G'_2=G_2-N_{G_2}[y_1]$. Observe that $N_{G}[\{u_1,u_2\}]=N_{G_1}[\{x_1,x_2\}]\times N_{G_2}[y_1]$, therefore, $V(G-N_{G}[\{u_1,u_2\}])=(\cup_{x\in V(G'_1)}{}_{x}G_2)\cup (\cup_{y\in V(G'_2)}G_{1y})$. It is sufficient to show that all elements of $(\cup_{x\in V(G'_1)}{}_{x}G_2)\cup (\cup_{y\in V(G'_2)}G_{1y})$ are contained in a component. Since $girth(G_i)\geq 5$ for $i=1,2$, $V(G'_1)\neq\emptyset$ and $V(G'_2)\neq\emptyset$. Let $x_3\in V(G'_1)$ and $y_2\in V(G'_2)$. Observe that $V({}_{x_3}G_2)$ and $V(G_{1y})$, for any $y\in V(G'_2)$, are contained in a component; and $V(G_{1y_2})$ and $V({}_{x}G_2)$, for any $x\in V(G'_1)$, are contained in a component. This implies that all elements of $(\cup_{x\in V(G'_1)}{}_{x}G_2)\cup (\cup_{y\in V(G'_2)}G_{1y})$ are contained in the same component. Therefore $G-N_{G}[\{u_1,u_2\}]$ is connected and $N_{G}(\{u_1,u_2\})$ is a $1$-$extra$ $cut$ in $G$. Then~$\kappa_1(G)\leq |N_{G}(\{u_1,u_2\})|=2k_1k_2+2k_1-2$. Analogously, if $u_1=(x_1,y_1)$ and $u_2=(x_1,y_2)$, we also get that $N_{G}(\{u_1,u_2\})$ is a $1$-$extra$ $cut$ in $G$ and $\kappa_1(G)\leq |N_{G}(\{u_1,u_2\})|= 2k_1k_2+2k_2-2$. Therefore, $\kappa_1(G)\leq$ min$\{2k_1k_2+2k_1-2,2k_1k_2+2k_2-2\}$.

Now, it is sufficient to prove $\kappa_1(G)\geq$ min$\{2k_1k_2+2k_1-2,2k_1k_2+2k_2-2\}$. Let $S$ be a $\kappa_{1}$-$cut$ in $G$. Denote ${}_{x}G'_2={}_{x}G_2-{}_{x}S$ for $x\in V(G_1)$. We consider two cases in the following.

\noindent{\bf Case 1. }${}_{x}S\neq\emptyset$ for all $x\in V(G_1)$, or $S_y\neq \emptyset$ for all $y\in V(G_2)$.

Assume ${}_{x}S\neq\emptyset$ for all $x\in V(G_1)$. By $Lemma\,\ref{0}$ and $Lemma\,\ref{4}$, $|S|= \sum_{x\in V(G_1)}|{}_{x}S|\geq |V(G_1)|\kappa(G_2)\geq (3(k_1-1)+2)k_2\geq 2k_1k_2+2k_1-2$. Analogously, if $S_y\neq \emptyset$ for all $y\in V(G_2)$, then $|S|=\sum_{y\in V(G_2)}|S_y|\geq |V(G_2)|\kappa(G_1)\geq(3(k_2-1)+2)k_1\geq 2k_1k_2+2k_2-2$.

\noindent{\bf Case 2. }There exist a vertex $x'\in V(G_1)$ and a vertex $y'\in V(G_2)$ such that ${}_{x'}S=S_{y'}=\emptyset$.

By the assumption, ${}_{x'}S=S_{y'}=\emptyset$, we know that $V({}_{x'}G_2)$ and $V(G_{1y'})$ are contained in a component $H'$ of $G-S$. Let $H$ be another component of $G-S$. Denote $p_1(V(H))=\{x'_{1},\cdots,x'_{l}\}$ and ${}_{x}H={}_{x}G_2\cap H$ for $x\in V(G_1)$. We claim that if $C_i$ is a component of ${}_{x'_{i}}H$ for $i\in \{1,\cdots,l\}$, then there are at least $\kappa(G_1)|N_{{}_{x'_{i}}G_2}[V(C_i)]|$ internally disjoint paths between $V(C_i)$ and $V({}_{x'}G'_{2})(=V({}_{x'}G_{2}))$ in $G-{}_{x'}S-{}_{x'_{i}}S$. Let $C_{i'}$ be a component of ${}_{x'_{i'}}H$ for some $i'\in \{1,\cdots,l\}$. Now, we introduce some general constructions of internally disjoint paths between $V(C_{i'})$ and $V({}_{x'}G'_{2})$ in $G-{}_{x'}S-{}_{x'_{i'}}S$. Let $u\in V(C_{i'})$. By $Theorem\,\ref{1}$, there are at least $\kappa(G_1)$ internally disjoint paths between $u$ and $(x',p_2(u))$ in $G_{1p_2(u)}$, therefore we have $\kappa(G_1)|V(C_{i'})|$ internally disjoint paths between $V(C_{i'})$ and $\{x'\}\times p_2(V(C_{i'}))$ in $G-{}_{x'}S-{}_{x'_{i'}}S$. To construct the $\kappa(G_1)|N_{{}_{x'_{i'}}G_2}(V(C_{i'}))|$ remaining paths, we consider internally disjoint paths between $v$ and $(x',p_2(v))$ in $G_{1p_2(v)}$ for $v\in N_{{}_{x'_{i'}}G_2}(V(C_{i'}))$. Let $v'\in N_{{}_{x'_{i'}}G_2}(V(C_{i'}))$. By $Theorem\,\ref{1}$, there are at least $\kappa(G_1)$ internally disjoint paths between $v'$ and $(x',p_2(v'))$ in $G_{1p_2(v')}$, we denote them by $P_1,P_2,\cdots,P_{\kappa(G_1)}$. Let $P'_s=P_s-v'$ for $s\in\{1,2,\cdots,\kappa(G_1)\}$. Since $v'\in N_{{}_{x'_{i'}}G_2}(V(C_{i'}))$, $v'$ has a neighbor denoted by $u'$ in $C_{i'}$. Observe that $u'P'_s$, for $s\in\{1,2,\cdots,\kappa(G_1)\}$, is a path between $u'$ and $(x',p_2(v'))$ in $G-{}_{x'}S-{}_{x'_{i'}}S$. This implies that there are $\kappa(G_1)|N_{{}_{x'_{i'}}G_2}(V(C_{i'}))|$ internally disjoint paths between $V(C_{i'})$ and $\{x'\}\times p_2(N_{{}_{x'_{i'}}G_2}(V(C_{i'})))$ in $G-{}_{x'}S-{}_{x'_{i'}}S$. Therefore, we have $\kappa(G_1)|N_{{}_{x'_{i'}}G_2}[V(C_{i'})]|$ internally disjoint paths between $V(C_{i'})$ and $V({}_{x'}G'_{2})$ in $G-{}_{x'}S-{}_{x'_{i'}}S$. Since $G-S$ is not connected, by $Theorem\,\ref{1}$, we know that $|S|\geq \kappa(G_1)|N_{{}_{x'_{i'}}G_2}[V(C_{i'})]|+ |{}_{x'_{i'}}S|$. It makes sense to consider subcases in the following.

\noindent{\bf Subcase 2.1. }$|V(C_{i'})|\geq 3(k_2-1)+2$ for some $i'\in\{1,\cdots,l\}$.

In this case, $|S|\geq \kappa(G_1)|N_{{}_{x'_{i'}}G_2}[V(C_{i'})]|+ |{}_{x'_{i'}}S|>\kappa(G_1)|V(C_{i'})|\geq (3(k_2-1)+2)k_1\geq 2k_1k_2+2k_2-2$.

\noindent{\bf Subcase 2.2. }$2\leq |V(C_{i'})|\leq 3(k_2-1)+1$ for some $i'\in\{1,\cdots,l\}$.

If $3\leq |V(C_{i'})|\leq 3(k_2-1)+1$, since $girth(G_2)\geq 5$, $|{}_{x'_{i'}}S|\geq |N_{{}_{x'_{i'}}G_{2}}(V(C_{i'}))|\geq |N_{{}_{x'_{i'}}G_{2}}[\{u_1,u_2, u_{3}\}]|-|V(C_{i'})|\geq 3k_2-1- |V(C_{i'})|$, where $C_{i'}[\{u_1,u_2, u_{3}\}]$ is connected.  Therefore $|S|\geq \kappa(G_1)|N_{{}_{x'_{i'}}G_{2}}[V(C_{i'})]|+|{}_{x'_{i'}}S|\geq (3k_2-1)k_1+1> 2k_1k_2+2k_2-2$.

If $|V(C_{i'})|=2$, since $girth(G_2)\geq 5$, $|{}_{x'_{i'}}S|\geq |N_{{}_{x'_{i'}}G_{2}}(V(C_{i'}))|= |N_{{}_{x'_{i'}}G_{2}}(\{u_1,u_2\})|=2k_2-2$, where $u_1u_{2}$ is an edge of $C_{i'}$. Therefore $|S|\geq \kappa(G_1)|N_{{}_{x'_{i'}}G_{2}}[V(C_i)]|+|{}_{x'_{i'}}S|\geq (2k_2-2+2)k_1 + 2k_2-2=2k_1k_2+2k_2-2$.

\noindent{\bf Subcase 2.3. }$|V(C_{i'})|= 1$ for some $i'\in\{1,\cdots,l\}$

 Since $S$ is a $\kappa_1$-$cut$ in $G$, $p_1(N_{H}(V(C_{i'})))\cap (p_1(V(H))\backslash \{x'_{i'}\})\neq \emptyset$. This implies $|p_1(V(H))|\geq 2$. Subcase 2.1 and 2.2 imply that if $C_i$ is a component of ${}_{x'_{i}}H$ for $i\in \{1,\cdots,l\}$ and $|V(C_{i})|\geq 2$, then $|S|\geq 2k_1k_2+2k_2-2$, hence, we suppose that any component of ${}_{x'_{i}}H$ is an isolated vertex for $i\in \{1,\cdots,l\}$ in the remaining proof. Since $S$ is a $\kappa_1$-$cut$ in $G$, $N_{G}(V(H))= S$. It makes sense to consider a lower bound of the order of $N_{G}(V(H))$.  Denote $D_i=N_{{}_{x'_{i}}G_2}(V({}_{x'_i}H))$ for $i\in\{1,\cdots,l\}$. Since $G_2$ is $k_2$-regular and maximally connected, $|D_i|\geq k_2$. Observe that $|N_{G}(V(H))|\geq \sum_{x'_{i}\in p_1(V(H))}|D_i|+|\cup^{l}_{i=1}(N_{G_1}(x'_i)\backslash p_1(V(H))\times p_2(N_{{}_{x'_i}G_2}[V({}_{x'_i}H)])|$. And we divide subcase 2.3 into the following subcases.

\noindent{\bf Subcase 2.3.1. }$|p_1(V(H))|\geq 3(k_1-1)+2$.

In this case, $|S|\geq \sum_{x'_{i}\in p_1(V(H))}|D_i|+|\cup^{l}_{i=1}(N_{G_1}(x'_i)\backslash p_1(V(H))\times p_2(N_{{}_{x'_i}G_2}[V({}_{x'_i}H)])|>\sum_{x'_{i}\in p_1(V(H))}|D_i|\geq |p_1(V(H))|\kappa(G_2)\geq (3(k_1-1)+2)k_2\geq 2k_1k_2+2k_1-2$.

\noindent{\bf Subcase 2.3.2. }$2\leq |p_1(V(H))|\leq 3(k_1-1)+1$.

Clearly, $H$ is connected, therefore $G_1[p_1(V(H))]$ is connected.

If $3\leq |p_1(V(H))|\leq 3(k_1-1)+1$, assume, without loss of generality, $G_1[\{x'_{1},x'_{2}, x'_{3}\}]$ is connected. Let $(x'_i,y'_i)\in V(H)$ for $i=1,2,3$. Since $girth(G_1)\geq 5$, $|N_{G_1}(p_1(V(H)))|\geq |N_{G_1}[\{x'_{1},x'_{2},x'_{3}\}]|-|p_1(V(H))|=3k_1-1-l$. It implies that there are at least $3k_1-1+2-l$ neighbors of $V(H)$ in $\cup^{3}_{i=1}(N_{G_1}[x'_i]\backslash p_1(V(H)))\times \{y'_i\}$. These neighbors are denoted by $\Rmnum{1}$-type. Observe that $\cup^{l}_{i=1}(N_{G_1}[x'_i]\backslash p_1(V(H)))\times p_2(D_i)\subseteqq N_G(V(H))$. And $\sum^{3}_{i=1}|(N_{G_1}[x'_i]\backslash p_1(V(H)))\times p_2(D_i)|\geq (3k_1-1-l)k_2$. Therefore we have at least $(3k_1-1-l)k_2$ neighbors of $V(H)$ contained in $\cup^{3}_{i=1}(N_{G_1}[x'_i]\backslash p_1(V(H)))\times p_2(D_i)$. These neighbors are denoted by $\Rmnum{2}$-type. Since $girth(G_1)\geq 5$, $\Rmnum{1}$-type and $\Rmnum{2} $-type do not have the same element. Therefore there are at least $3k_1-1-l+(3k_1-1-l)k_2$ neighbors of $V(H)$ in $\cup^{l}_{i=1}(N_{G_1}(x'_i)\backslash p_1(V(H))\times p_2(N_{{}_{x'_i}G_2}[V({}_{x'_i}H)])$. Then $|S|\geq |N_{G}(V(H))|\geq \sum_{x'_{i}\in p_1(V(H))}|D_i|+|\cup^{l}_{i=1}(N_{G_1}(x'_i)\backslash p_1(V(H))\times p_2(N_{{}_{x'_i}G_2}[V({}_{x'_i}H)])|\geq lk_2+3k_1-1-l+(3k_1-1-l)k_2\geq (3k_1-1)k_2+1> 2k_1k_2+2k_1-2$ holds.

If $|p_1(V(H))|=2$, assume, without loss of generality, $x'_{1}x'_{2}$ is an edge in $G_1[p_1(V(H))]$.
Since $girth(G_1)\geq 5$, $|N_{G_1}(p_1(V(H)))|= |N_{G_1}(\{x'_{1},x'_{2}\})|=2k_1-2$. This implies that there are $2(k_1-2)+2$ vertices of $\Rmnum{1}$-type and $(2(k_1-2)+2)k_2$ vertices of $\Rmnum{2}$-type. Therefore, there are at least $2k_1-2+(2k_1-2)k_2$ neighbors of $V(H)$ in $\cup^{l}_{i=1}(N_{G_1}(x'_i)\backslash p_1(V(H))\times p_2(N_{{}_{x'_i}G_2}[V({}_{x'_i}H)])$. Then $|S|\geq |N_{G}(V(H))|\geq \sum_{x'_{i}\in p_1(V(H))}|D_i|+|\cup^{l}_{i=1}(N_{G_1}(x'_i)\backslash p_1(V(H))\times p_2(N_{{}_{x'_i}G_2}[V({}_{x'_i}H)])|\geq 2k_2+ 2k_1-2+(2k_1-2)k_2= 2k_1k_2+2k_1-2$ holds.

According to the argument of subcase 2.3.1. and subcase 2.3.2., we can get that if any component of ${}_{x'_{i}}H$ is an isolated vertex for $i\in \{1,\cdots,l\}$, then $|S|\geq |N_G(V(H))|\geq 2k_1k_2+2k_1-2$. Therefore, if $C_i$ is a component of ${}_{x'_{i}}H$ for $i\in \{1,\cdots,l\}$, then $|S|\geq$ min$\{2k_1k_2+2k_1-2,2k_1k_2+2k_2-2\}$ holds for all cases of $|V(C_i)|$. Combining with case 1, $\kappa_1(G)\geq$ min$\{2k_1k_2+2k_1-2,2k_1k_2+2k_2-2\}$. The proof is thus complete. $\Box$

\begin{theorem}\label{7}
Let $G_i$ be a $k_i\ (\geq 2)$-regular and maximally connected graph with $girth(G_i)\geq 6$ for $i=1,2$. Then $\kappa_{2}(G_1\boxtimes G_2)=$ min$\{3k_1k_2+3k_1-k_2-4,3k_1k_2-k_1+3k_2-4\}$.
\end{theorem}
\noindent{\bf Proof.} Denote $G=G_1\boxtimes G_2$. Let $u_1u_2u_3$ be a path in $G$, where $u_1=(x_1,y_1)$, $u_2=(x_2,y_1)$ and $u_3=(x_3,y_1)$. It is routine to verify that $N_{G}(\{u_1,u_2, u_3\})$ is a $2$-$extra$ $cut$ in $G$. Clearly, $|G-N_{G}[\{u_1,u_2, u_3\}]|> 3$, therefore if we can prove that $G-N_{G}[\{u_1,u_2, u_3\}]$ is connected, then $N_{G}(\{u_1,u_2, u_3\})$ is a $2$-$extra$ $cut$ in $G$. Denote $G'_1=G_1-N_{G_1}[\{x_1,x_2,x_3\}]$ and $G'_2=G_2-N_{G_2}[y_1]$. Observe that $N_{G}[\{u_1,u_2,u_3\}]=N_{G_1}[\{x_1,x_2,x_3\}]\times N_{G_2}[y_1]$, therefore, $V(G-N_{G}[\{u_1,u_2,u_3\}])=(\cup_{x\in V(G'_1)}{}_{x}G_2)\cup (\cup_{y\in V(G'_2)}G_{1y})$. It is sufficient to show that all elements of $(\cup_{x\in V(G'_1)}{}_{x}G_2)\cup (\cup_{y\in V(G'_2)}G_{1y})$ are contained in a component. Since $girth(G_i)\geq 6$ for $i=1,2$, $V(G'_1)\neq\emptyset$ and $V(G'_2)\neq\emptyset$. Let $x_4\in V(G'_1)$ and $y_2\in V(G'_2)$. Observe that $V({}_{x_4}G_2)$ and $V(G_{1y})$ for any $y\in V(G'_2)$ are contained in a component; and $V(G_{1y_2})$ and $V({}_{x}G_2)$ for any $x\in V(G'_1)$ are contained in a component. This implies that all elements of $(\cup_{x\in V(G'_1)}{}_{x}G_2)\cup (\cup_{y\in V(G'_2)}G_{1y})$ are contained in the same component. Therefore $G-N_{G}[\{u_1,u_2,u_3\}]$ is connected and  $N_{G}(\{u_1,u_2,u_3\})$ is a $2$-$extra$ $cut$ in $G$. Then $\kappa_2(G)\leq |N_{G}(\{u_1,u_2,u_3\})| =3k_1k_2+3k_1-k_2-4$. Analogously, if $u_1=(x_1,y_1)$, $u_2=(x_1,y_2)$ and $u_3=(x_1,y_3)$, we also get that $N_{G}(\{u_1,u_2,u_3\})$ is a $2$-$extra$ $cut$ in $G$ and $\kappa_2(G)\leq |N_{G}(\{u_1,u_2,u_3\})|= 3k_1k_2-k_1+3k_2-4$. Therefore, $\kappa_2(G)\leq$ min$\{3k_1k_2+3k_1-k_2-4,3k_1k_2-k_1+3k_2-4\}$.

Now, it is sufficient to prove $\kappa_2(G)\geq$ min$\{3k_1k_2+3k_1-k_2-4,3k_1k_2-k_1+3k_2-4\}$. Observe that min$\{3k_1k_2+3k_1-k_2-4,3k_1k_2-k_1+3k_2-4\}=3k_1k_2+3k_1-k_2-4$ if and only if $k_1\leq k_2$, and min$\{3k_1k_2+3k_1-k_2-4,3k_1k_2-k_1+3k_2-4\}=3k_1k_2-k_1+3k_2-4$ if and only if $k_2\leq k_1$. Let $S$ be a $\kappa_{2}$-$cut$ in $G$. Therefore, we need to prove that $|S|\geq 3k_1k_2+3k_1-k_2-4$ when $k_1\leq k_2$ and $|S|\geq 3k_1k_2-k_1+3k_2-4$ when $k_2\leq k_1$. Denote ${}_{x}G'_2={}_{x}G_2-{}_{x}S$ for $x\in V(G_1)$. We consider two cases in the following.

\noindent{\bf Case 1. }${}_{x}S\neq\emptyset$ for all $x\in V(G_1)$, or $S_y\neq \emptyset$ for all $y\in V(G_2)$.

Assume ${}_{x}S\neq\emptyset$ for all $x\in V(G_1)$. By $Lemma\,\ref{0}$ and $Lemma\,\ref{4}$, If $k_1\leq k_2$,  $|S|= \sum_{x\in V(G_1)}|{}_{x}S|\geq |V(G_1)|\kappa(G_2)\geq (4(k_1-1)+2)k_2\geq 3k_1k_2+3k_1-k_2-4$; if $k_2\leq k_1$,  $|S|= \sum_{x\in V(G_1)}|{}_{x}S|\geq |V(G_1)|\kappa(G_2)\geq (4(k_1-1)+2)k_2\geq 3k_1k_2-k_1+3k_2-4$. Analogously, if $S_y\neq \emptyset$ for all $y\in V(G_2)$, when $k_1\leq k_2$, $|S|=\sum_{y\in V(G_2)}|S_y|\geq |V(G_2)|\kappa(G_1)\geq(4(k_2-1)+2)k_1\geq 3k_1k_2+3k_1-k_2-4$; when $k_2\leq k_1$, $|S|=\sum_{y\in V(G_2)}|S_y|\geq |V(G_2)|\kappa(G_1)\geq(4(k_2-1)+2)k_1\geq 3k_1k_2-k_1+3k_2-4$.

\noindent{\bf Case 2. }There exist a vertex $x'\in V(G_1)$ and a vertex $y'\in V(G_2)$ such that ${}_{x'}S=S_{y'}=\emptyset$.

By the assumption, ${}_{x'}S=S_{y'}=\emptyset$, we know that $V({}_{x'}G_2)$ and $V(G_{1y'})$ are contained in a component $H'$ of $G-S$. Let $H$ be another component of $G-S$. Denote $p_1(V(H))=\{x'_{1},\cdots,x'_{l}\}$ and ${}_{x}H={}_{x}G_2\cap H$ for $x\in V(G_1)$. Doing the same argument as $Theorem\,\ref{6}$, we can get that if $C_i$ is a component of ${}_{x'_{i}}H$ for $i\in \{1,\cdots,l\}$, then there are at least $\kappa(G_1)|N_{{}_{x'_{i}}G_2}[V(C_i)]|$ internally disjoint paths between $V(C_i)$ and $V({}_{x'}G'_{2})(=V({}_{x'}G_{2}))$ in $G-{}_{x'}S-{}_{x'_{i}}S$. And then we consider subcases in the following.

\noindent{\bf Subcase 2.1. }$|V(C_{i'})|\geq 4(k_2-1)+2$ for some $i'\in\{1,\cdots,l\}$.

In this case, $|S|\geq \kappa(G_1)|N_{{}_{x'_{i'}}G_2}[V(C_{i'})]|+ |{}_{x'_{i'}}S|>\kappa(G_1)|V(C_{i'})| \geq (4(k_2-1)+2)k_1$. Therefore, $|S|\geq 3k_1k_2+3k_1-k_2-4$ when $k_1\leq k_2$, and $|S|\geq 3k_1k_2-k_1+3k_2-4$ when $k_2\leq k_1$.

\noindent{\bf Subcase 2.2. }$3\leq|V(C_{i'})|\leq 4(k_2-1)+1$ for some $i'\in\{1,\cdots,l\}$.

If $4\leq|V(C_{i'})|\leq 4(k_2-1)+1$, then
$|{}_{x'_{i'}}S|\geq |N_{{}_{x'_{i'}}G_{2}}(V(C_{i'}))|\geq |N_{{}_{x'_{i'}}G_{2}}[\{u_1,u_2, u_{3},u_4\}]|-|V(C_{i'})|$, where $C_{i'}[\{u_1,u_2,u_{3},u_4\}]$ is connected. Since $girth(G_2)\geq 6$, it is easy to verify that $|N_{{}_{x'_{i'}}G_{2}}[V(C_{i'})]|\geq |N_{{}_{x'_{i'}}G_{2}}[\{u_1,u_2, u_{3},u_4\}]|\geq 4k_2-2$. Then $|S|\geq \kappa(G_1)|N_{{}_{x'_{i'}}G_{2}}[V(C_{i'})]|+|{}_{x'_{i'}}S|\geq (4k_2-2)k_1+1$. Therefore, $|S|> 3k_1k_2+3k_1-k_2-4$ when $k_1\leq k_2$, and $|S|> 3k_1k_2-k_1+3k_2-4$ when $k_2\leq k_1$.

If $|V(C_{i'})|=3$, since $girth(G_2)\geq 6$, $|{}_{x'_{i'}}S|\geq|N_{{}_{x'_{i'}}G_{2}}(V(C_{i'}))|= |N_{{}_{x'_{i'}}G_{2}}(\{u_1,u_2,u_3\})|=3k_2-4$, where $C_{i'}[\{u_1,u_{2},u_3\}]$ is connected. This implies that $|S|\geq \kappa(G_1)(|N_{{}_{x'_{i'}}G_{2}}[V(C_{i'})]|+|{}_{x'_{i'}}S|\geq (3k_2-4+3)k_1 + 3k_2-4=3k_1k_2-k_1+3k_2-4$. Therefore, $|S|\geq 3k_1k_2+3k_1-k_2-4$ when $k_1\leq k_2$, and $|S|\geq 3k_1k_2-k_1+3k_2-4$ when $k_2\leq k_1$.

\noindent{\bf Subcase 2.3. }$|V(C_{i'})|\leq 2$ for some $i'\in\{1,\cdots,l\}$.

Since $S$ is a $\kappa_2$-$cut$ in $G$, $p_1(N_{H}(V(C_{i'})))\cap (p_1(V(H))\backslash \{x'_{i'}\})\neq \emptyset$. This implies $|p_1(V(H))|\geq 2$. Subcase 2.1 and 2.2 imply that if $C_i$ is a component of ${}_{x'_{i}}H$ for $i\in \{1,\cdots,l\}$ and $|V(C_{i})|\geq 3$, then $|S|\geq 3k_1k_2+3k_1-k_2-4$ when $k_1\leq k_2$ and $|S|\geq 3k_1k_2-k_1+3k_2-4$ when $k_2\leq k_1$, hence, we suppose that any component of ${}_{x'_{i}}H$ contains at most 2 vertices for $i\in \{1,\cdots,l\}$ in the remaining proof. Since $S$ is a $\kappa_2$-$cut$ in $G$, $N_{G}(V(H))= S$. It makes sense to consider a lower bound of the order of $N_{G}(V(H))$.  Denote $D_i=N_{{}_{x'_{i}}G_2}(V({}_{x'_i}H))$ for $i\in\{1,\cdots,l\}$. Since $G_2$ is $k_2$-regular and maximally connected, $|D_i|\geq k_2$. Observe that $|N_{G}(V(H))|\geq \sum_{x'_{i}\in p_1(V(H))}|D_i|+|\cup^{l}_{i=1}(N_{G_1}(x'_i)\backslash p_1(V(H))\times p_2(N_{{}_{x'_i}G_2}[V({}_{x'_i}H)])|$. And we divide subcase 2.3 into the following subcases.

\noindent{\bf Subcase 2.3.1. }$|p_1(V(H))|\geq 4(k_1-1)+2$.

In this case, $|S|\geq \sum_{x'_{i}\in p_1(V(H))}|D_i|+|\cup^{l}_{i=1}(N_{G_1}(x'_i)\backslash p_1(V(H))\times p_2(N_{{}_{x'_i}G_2}[V({}_{x'_i}H)])|>\sum_{x'_{i}\in p_1(V(H))}|D_i|\geq |p_1(V(H))|\kappa(G_2)\geq (4(k_1-1)+2)k_2$. Therefore, $|S|\geq 3k_1k_2+3k_1-k_2-4$ when $k_1\leq k_2$, and $|S|\geq 3k_1k_2-k_1+3k_2-4$ when $k_2\leq k_1$.

\noindent{\bf Subcase 2.3.2. }$3\leq |p_1(V(H))|\leq 4(k_1-1)+1$.

Clearly, $H$ is connected, therefore $G_1[p_1(V(H))]$ is connected.

If $4\leq |p_1(V(H))|\leq 4(k_1-1)+1$, assume, without loss of generality, $G_1[\{x'_{1},x'_{2},x'_{3},x'_{4}\}]$ is connected. Let $(x'_i,y'_i)\in V(H)$ for $i=1,2,3,4$. Since $girth(G_1)\geq 6$, $|N_{G_1}(p_1(V(H)))|\geq |N_{G_1}[\{x'_{1},x'_{2},x'_{3},x'_{4}\}]|-|p_1(V(H))|\geq 4k_1-2-l$. It implies that there are at least $4k_1-2-l$ neighbors of $V(H)$ in $\cup^{4}_{i=1}(N_{G_1}[x'_i]\backslash p_1(V(H)))\times \{y'_i\}$. These neighbors are denoted by $\Rmnum{1}$-type. Observe that $\cup^{l}_{i=1}(N_{G_1}[x'_i]\backslash p_1(V(H)))\times p_2(D_i)\subseteqq N_G(V(H))$. And $\sum^{4}_{i=1}|(N_{G_1}[x'_i]\backslash p_1(V(H)))\times p_2(D_i)|\geq (4k_1-2-l)k_2$. Therefore we have at least $(4k_1-2-l)k_2$ neighbors of $V(H)$ contained in $\cup^{4}_{i=1}(N_{G_1}[x'_i]\backslash p_1(V(H)))\times p_2(D_i)$. These neighbors are denoted by $\Rmnum{2}$-type. Since $girth(G_1)\geq 6$, $\Rmnum{1}$-type and $\Rmnum{2} $-type do not have the same element. Therefore there are at least $4k_1-2-l+(4k_1-2-l)k_2$ neighbors of $V(H)$ in $\cup^{l}_{i=1}(N_{G_1}(x'_i)\backslash p_1(V(H))\times p_2(N_{{}_{x'_i}G_2}[V({}_{x'_i}H)])$. Then $|S|\geq |N_{G}(V(H))|\geq \sum_{x'_{i}\in p_1(V(H))}|D_i|+|\cup^{l}_{i=1}(N_{G_1}(x'_i)\backslash p_1(V(H))\times p_2(N_{{}_{x'_i}G_2}[V({}_{x'_i}H)])|\geq lk_2+4k_1-2-l+(4k_1-2-l)k_2\geq (4k_1-2)k_2+1$. Therefore, $|S|> 3k_1k_2+3k_1-k_2-4$ when $k_1\leq k_2$, and $|S|\geq 3k_1k_2-k_1+3k_2-4$ when $k_2\leq k_1$.

If $|p_1(V(H))|=3$, assume, without loss of generality, $G_1[\{x'_{1},x'_{2},x'_{3}\}]$ is connected.
Since $girth(G_1)\geq 6$, $|N_{G_1}(p_1(V(H)))|= |N_{G_1}(\{x'_{1},x'_{2},x'_{3}\})|=3k_1-4$. This implies that there are $3k_1-4$ vertices of $\Rmnum{1}$-type and $(3k_1-4)k_2$ vertices of $\Rmnum{2}$-type. Therefore, there are at least $3k_1-4+(3k_1-4)k_2$ neighbors of $V(H)$ in $\cup^{l}_{i=1}(N_{G_1}(x'_i)\backslash p_1(V(H))\times p_2(N_{{}_{x'_i}G_2}[V({}_{x'_i}H)])$. Then $|S|\geq |N_{G}(V(H))|\geq \sum_{x'_{i}\in p_1(V(H))}|D_i|+|\cup^{l}_{i=1}(N_{G_1}(x'_i)\backslash p_1(V(H))\times p_2(N_{{}_{x'_i}G_2}[V({}_{x'_i}H)])|\geq 3k_2+ 3k_1-4+(3k_1-4)k_2= 3k_1k_2+3k_1-k_2-4$. Therefore, $|S|\geq 3k_1k_2+3k_1-k_2-4$ when $k_1\leq k_2$, and $|S|\geq 3k_1k_2-k_1+3k_2-4$ when $k_2\leq k_1$.

\begin{figure}[htbp]
\centering
\includegraphics[width=7cm]{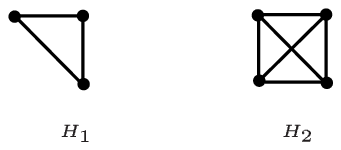}\\
\caption{Types of $H$ when $|p_1(V(H))|=2$ and $|p_2(V(H))|=2$.}
\label{10}
\end{figure}

According to the argument of subcase 2.3.1. and subcase 2.3.2., we can get that if any component of ${}_{x'_{i}}H$ contains at most 2 vertices for $i\in \{1,\cdots,l\}$ and $|p_1(V(H))|\geq 3$, then $|S|\geq 3k_1k_2+3k_1-k_2-4$ when $k_1\leq k_2$, and $|S|\geq 3k_1k_2-k_1+3k_2-4$ when $k_2\leq k_1$. Therefore, it remains to show that if any component of ${}_{x'_{i}}H$ contains at most 2 vertices for $i\in \{1,\cdots,l\}$ and $|p_1(V(H))|\leq 2$, then $|S|\geq 3k_1k_2+3k_1-k_2-4$ when $k_1\leq k_2$, and $|S|\geq 3k_1k_2-k_1+3k_2-4$ when $k_2\leq k_1$. Denote $H_y=G_{1y}\cap H$ for $y\in V(G_2)$, $p_2(V(H))=\{y'_1,\cdots,y'_h\}$ and $T_j=N_{G_{1y'_j}}(V(H_{y'_j}))$ for $j\in\{1,\cdots,h\}$. Since $|p_1(V(H)))|\leq 2$, $|V(H_{y'_j})|\leq 2$ for $j\in\{1,\cdots,h\}$. Doing similar arguments as subcase 2.3.1. and subcase 2.3.2., we can get that if $|p_2(V(H))|\geq 3$, then $|S|\geq 3k_1k_2+3k_1-k_2-4$ when $k_1\leq k_2$, and $|S|\geq 3k_1k_2-k_1+3k_2-4$ when $k_2\leq k_1$. Since $S$ is a $\kappa_2$-$cut$ in $G$, we only need to show in the remaining proof that if $|p_1(V(H))|=2$ and $|p_2(V(H))|=2$ then $|S|\geq 3k_1k_2+3k_1-k_2-4$ when $k_1\leq k_2$, and $|S|\geq 3k_1k_2-k_1+3k_2-4$ when $k_2\leq k_1$. And there are two types of $H$ as shown in Figure\ref{10}.

Observe that $|N_{G}(V(H_1))|=3k_1k_2+k_1+k_2-4\geq 3k_1k_2+3k_1-k_2-4$ when $k_1\leq k_2$, and $|N_{G}(V(H_1))|=3k_1k_2+k_1+k_2-4\geq 3k_1k_2-k_1+3k_2-4$ when $k_2\leq k_1$; $ |N_{G}(V(H_2))|=4k_1k_2-4\geq 3k_1k_2+3k_1-k_2-4$ when $k_1\leq k_2$, and $|N_{G}(V(H_1))|=4k_1k_2-4\geq 3k_1k_2-k_1+3k_2-4$ when $k_2\leq k_1$. Therefore, if $C_i$ is a component of ${}_{x'_{i}}H$ for $i\in \{1,\cdots,l\}$, then $|S|\geq$ min$\{3k_1k_2+3k_1-k_2-4,3k_1k_2-k_1+3k_2-4\}$ holds for all case of $|V(C_i)|$. Combining with case 1, Thus $\kappa_2(G)\geq$ min$\{3k_1k_2+3k_1-k_2-4,3k_1k_2-k_1+3k_2-4\}$, the proof is complete. $\Box$

\begin{theorem}\label{8}
Let $G_i$ be a $k_i\ (\geq 2)$-regular and maximally connected graph with $girth(G_i)\geq 7$ for $i=1,2$. Then $\kappa_{3}(G_1\boxtimes G_2)=$ min$\{4k_1k_2+4k_1-2k_2-6,4k_1k_2-2k_1+4k_2-6,4k_1k_2-4\}$.	
\end{theorem}
\noindent{\bf Proof.} Denote $G=G_1\boxtimes G_2$. Let $x_1x_2x_3x_4$ be a path in $G_1$ and $y_1y_2y_3y_4$ be a path in $G_2$. Denote $A_1=\{x_1,x_2,x_3,x_{4}\}\times\{y_1\}$, $A_2=\{x_1\}\times\{y_1,y_2,y_3,y_{4}\}$. It is routine to verify that $N_{G}(A_i)$ is a $3$-$extra$ $cut$ in $G$ for $i=1,2$. We only need to show that $N_G(A_1)$ is a $3$-$extra$ $cut$ in $G$ since the proof of $N_G(A_2)$ is similar to it. Clearly, $|G-N_{G}[A_1]|> 4$, therefore if we can prove that $G-N_{G}[A_1]$ is connected, then $N_{G}(A_1)$ is a $3$-$extra$ $cut$ in $G$. Denote $G'_1=G_1-N_{G_1}[\{x_1,x_2,x_3,x_4\}]$ and $G'_2=G_2-N_{G_2}[y_1]$. Observe that $N_{G}[A_1]=N_{G_1}[\{x_1,x_2,x_3,x_4\}]\times N_{G_2}[y_1]$, therefore, $V(G-N_{G}[A_1])=(\cup_{x\in V(G'_1)}{}_{x}G_2)\cup (\cup_{y\in V(G'_2)}G_{1y})$. It is sufficient to show that all elements of $(\cup_{x\in V(G'_1)}{}_{x}G_2)\cup (\cup_{y\in V(G'_2)}G_{1y})$ are contained in a component. Since $girth(G_i)\geq 7$ for $i=1,2$, $V(G'_1)\neq\emptyset$ and $V(G'_2)\neq\emptyset$. Let $x_5\in V(G'_1)$ and $y_2\in V(G'_2)$. Observe that $V({}_{x_5}G_2)$ and $V(G_{1y})$ for any $y\in V(G'_2)$ are contained in a component; and $V(G_{1y_2})$ and $V({}_{x}G_2)$ for any $x\in V(G'_1)$ are contained in a component. This implies that all elements of $(\cup_{x\in V(G'_1)}{}_{x}G_2)\cup (\cup_{y\in V(G'_2)}G_{1y})$ are contained in the same component. Therefore $G-N_{G}[A_1]$ is connected and $N_{G}(A_1)$ is a $3$-$extra$ $cut$ in $G$. Then $\kappa_3(G)\leq |N_{G}(A_1)| =4k_1k_2+4k_1-2k_2-6$. Analogously, we can also get that $N_{G}(A_2)$ is a $3$-$extra$ $cut$. If $k_2\geq 2k_1-1$ or $k_1\geq 2k_2-1$,   $4k_1k_2-4\geq$ min$\{4k_1k_2+4k_1-2k_2-6,4k_1k_2-2k_1+4k_2-6\}$, therefore, $\kappa_3(G)\leq$ min$\{4k_1k_2+4k_1-2k_2-6,4k_1k_2-2k_1+4k_2-6\}$. If $k_1\leq k_2\leq 2k_1-1$ or $k_2\leq k_1\leq 2k_2-1$, let $A_3=\{x_1,x_2\}\times\{y_1,y_{2}\}$, using the same argument as above, we can also verify that $N_{G}(A_3)$ is a $3$-$extra$ $cut$ in $G$. Therefore, $\kappa_3(G)\leq$ min$\{4k_1k_2+4k_1-2k_2-6,4k_1k_2-2k_1+4k_2-6,4k_1k_2-4\}$.

Now, it is sufficient to prove $\kappa_2(G)\geq$ min$\{4k_1k_2+4k_1-2k_2-6,4k_1k_2-2k_1+4k_2-6,4k_1k_2-4\}$. Observe that min$\{4k_1k_2+4k_1-2k_2-6,4k_1k_2-2k_1+4k_2-6,4k_1k_2-4\}=4k_1k_2-4$ if and only if $k_1\leq k_2\leq 2k_1-1$ or $k_2\leq k_1\leq 2k_2-1$, min$\{4k_1k_2+4k_1-2k_2-6,4k_1k_2-2k_1+4k_2-6,4k_1k_2-4\}=4k_1k_2+4k_1-2k_2-6$ if and only if $k_2\geq 2k_1-1$ and min$\{4k_1k_2+4k_1-2k_2-6,4k_1k_2-2k_1+4k_2-6,4k_1k_2-4\}=4k_1k_2-2k_1+4k_2-6$ if and only if $k_1\geq 2k_2-1$. Let $S$ be a $\kappa_{3}$-$cut$ in $G$. Therefore, we need to prove that $|S|\geq 4k_1k_2-4$ when $k_1\leq k_2\leq 2k_1-1$ or $k_2\leq k_1\leq 2k_2-1$, $|S|\geq 4k_1k_2+4k_1-2k_2-6$ when $k_2\geq 2k_1-1$, and $|S|\geq 4k_1k_2-2k_1+4k_2-6$ when $k_1\geq 2k_2-1$. Denote ${}_{x}G'_2={}_{x}G_2-{}_{x}S$ for $x\in V(G_1)$. We consider two cases in the following.

\noindent{\bf Case 1. } ${}_{x}S\neq\emptyset$ for all $x\in V(G_1)$, or $S_y\neq \emptyset$ for all $y\in V(G_2)$.

Assume ${}_{x}S\neq\emptyset$ for all $x\in V(G_1)$. By $Lemma\,\ref{0}$ and $Lemma\,\ref{4}$, if $k_1\leq k_2\leq 2k_1-1$ or $k_2\leq k_1\leq 2k_2-1$, $|S|= \sum_{x\in V(G_1)}|{}_{x}S|\geq |V(G_1)|\kappa(G_2)\geq (5(k_1-1)+2)k_2>4k_1k_2-4$; if $k_2\geq 2k_1-1$, $|S|= \sum_{x\in V(G_1)}|{}_{x}S|\geq |V(G_1)|\kappa(G_2)\geq(5(k_1-1)+2)k_2>4k_1k_2+4k_1-2k_2-6$; if $k_1\geq 2k_2-1$, $|S|= \sum_{x\in V(G_1)}|{}_{x}S|\geq |V(G_1)|\kappa(G_2)\geq(5(k_1-1)+2)k_2>4k_1k_2-2k_1+4k_2-6$. Analogously, if $S_y\neq \emptyset$ for all $y\in V(G_2)$, when $k_1\leq k_2\leq 2k_1-1$ or $k_2\leq k_1\leq 2k_2-1$, $|S|=\sum_{y\in V(G_2)}|S_y|\geq |V(G_2)|\kappa(G_1)\geq(5(k_2-1)+2)k_1>4k_1k_2-4$; when $k_2\geq 2k_1-1$, $|S|=\sum_{y\in V(G_2)}|S_y|\geq |V(G_2)|\kappa(G_1)\geq(5(k_2-1)+2)k_1>4k_1k_2+4k_1-2k_2-6$; when $k_1\geq 2k_2-1$, $|S|=\sum_{y\in V(G_2)}|S_y|\geq |V(G_2)|\kappa(G_1)\geq(5(k_2-1)+2)k_1>4k_1k_2-2k_1+4k_2-6$.

\noindent{\bf Case 2. } There exist a vertex $x'\in V(G_1)$ and a vertex $y'\in V(G_2)$ such that ${}_{x'}S=S_{y'}=\emptyset$.

By the assumption, ${}_{x'}S=S_{y'}=\emptyset$, we know that $V({}_{x'}G_2)$ and $V(G_{1y'})$ are contained in a component $H'$ of $G-S$. Let $H$ be another component of $G-S$. Denote $p_1(V(H))=\{x'_{1},\cdots,x'_{k}\}$ and ${}_{x}H={}_{x}G_2\cap H$ for $x\in V(G_1)$. Doing the same argument as $Theorem\,\ref{6}$, we can get that if $C_i$ is a component of ${}_{x'_{i}}H$ for $i\in \{1,\cdots,l\}$, then there are at least $\kappa(G_1)|N_{{}_{x'_{i}}G_2}[V(C_i)]|$ internally disjoint paths between $V(C_i)$ and $V({}_{x'}G'_{2})(=V({}_{x'}G_{2}))$ in $G-{}_{x'}S-{}_{x'_{i}}S$. And then we consider subcases in the following.

\noindent{\bf Subcase 2.1. }$|V(C_{i'})|\geq 5(k_2-1)+2$ for some $i'\in\{1,\cdots,l\}$.

In this case, $|S|\geq \kappa(G_1)|N_{{}_{x'_{i'}}G_2}[V(C_{i'})]|+ |{}_{x'_{i'}}S|> \kappa(G_1)|V(C_i)|\geq (5(k_2-1)+2)k_1$. Therefore, $|S|>4k_1k_2-4$ when $k_1\leq k_2\leq 2k_1-1$ or $k_2\leq k_1\leq 2k_2-1$, $|S|>4k_1k_2+4k_1-2k_2-6$ when $k_2\geq 2k_1-1$, and $|S|>4k_1k_2-2k_1+4k_2-6$ when $k_1\geq 2k_2-1$.

\noindent{\bf Subcase 2.2. }$4\leq|V(C_{i'})|\leq 5(k_2-1)+1$ for some $i'\in\{1,\cdots,l\}$.

If $5\leq|V(C_{i'})|\leq 5(k_2-1)+1$, then
$|N_{{}_{x'_{i'}}G_{2}}(V(C_{i'}))|\geq |N_{{}_{x'_{i'}}G_{2}}[\{u_1,u_2, u_{3},u_4,u_5\}]|-|V(C_{i'})|=5(k_2-1)+2-|V(C_{i'})|$, where $C_{i'}[\{u_1,u_2,u_{3},u_4,u_5\}]$ is connected. Since $girth(G_2)\geq 7$, it is easy to verify that $|N_{{}_{x'_{i'}}G_{2}}[V(C_{i'})]|\geq|N_{{}_{x'_{i'}}G_{2}}[\{u_1,u_2, u_{3},u_4,u_5\}]|\geq 5k_2-3$. Then $|S|\geq \kappa(G_1)|N_{{}_{x'_{i'}}G_{2}}[V(C_{i'})]|+|{}_{x'_{i'}}S|\geq (5k_2-3)k_1+1$. Therefore, $|S|>4k_1k_2-4$ when $k_1\leq k_2\leq 2k_1-1$ or $k_2\leq k_1\leq 2k_2-1$, $|S|>4k_1k_2+4k_1-2k_2-6$ when $k_2\geq 2k_1-1$, and $|S|>4k_1k_2-2k_1+4k_2-6$ when $k_1\geq 2k_2-1$.

If $|V(C_{i'})|=4$, then $|{}_{x'_{i'}}S|\geq|N_{{}_{x'_{i'}}G_{2}}(V(C_{i'}))|= |N_{{}_{x'_{i'}}G_{2}}(\{u_1,u_2,u_3,u_4\})|$, where $C_{i'}[\{u_1,u_{2},u_3,u_4\}]$ is connected. Since $girth(G_2)\geq 7$, it is easy to verify that $|N_{{}_{x'_{i'}}G_{2}}(\{u_1,u_2,u_3,u_4\})|\geq 4k_2-6$. Then $|S|\geq \kappa(G_1)|N_{{}_{x'_{i'}}G_{2}}[V(C_{i'})]|+|{}_{x'_{i'}}S|\geq (4k_2-6+4)k_1 + 4k_2-6=4k_1k_2-2k_1+4k_2-6$. Therefore, $|S|\geq 4k_1k_2-4$ when $k_1\leq k_2\leq 2k_1-1$ or $k_2\leq k_1\leq 2k_2-1$, $|S|> 4k_1k_2+4k_1-2k_2-6$ when $k_2\geq 2k_1-1$, and $|S|\geq4k_1k_2-2k_1+4k_2-6$ when $k_1\geq 2k_2-1$.

\noindent{\bf Subcase 2.3. }$|V(C_{i'})|\leq 3$ for some $i'\in\{1,\cdots,l\}$.

Since $S$ is a $\kappa_3$-$cut$ in $G$, $p_1(N_{H}(V(C_{i'})))\cap (p_1(V(H))\backslash \{x'_{i'}\})\neq \emptyset$. This implies $|p_1(V(H))|\geq 2$. Subcase 2.1 and 2.2 imply that if $C_i$ is a component of ${}_{x'_{i}}H$ for $i\in \{1,\cdots,l\}$ and $|V(C_{i})|\geq 4$, then $|S|\geq 4k_1k_2-4$ when $k_1\leq k_2\leq 2k_1-1$ or $k_2\leq k_1\leq 2k_2-1$, $|S|> 4k_1k_2+4k_1-2k_2-6$ when $k_2\geq 2k_1-1$, and $|S|\geq4k_1k_2-2k_1+4k_2-6$ when $k_1\geq 2k_2-1$, hence, we suppose that any component of ${}_{x'_{i}}H$ contains at most 3 vertices for $i\in \{1,\cdots,l\}$ in the remaining proof. Since $S$ is a $\kappa_3$-$cut$ in $G$, $N_{G}(V(H))= S$. It makes sense to consider a lower bound of the order of $N_{G}(V(H))$.  Denote $D_i=N_{{}_{x'_{i}}G_2}(V({}_{x'_i}H))$ for $i\in\{1,\cdots,l\}$. Since $G_2$ is $k_2$-regular and maximally connected, $|D_i|\geq k_2$. Observe that $|N_{G}(V(H))|\geq \sum_{x'_{i}\in p_1(V(H))}|D_i|+|\cup^{l}_{i=1}(N_{G_1}(x'_i)\backslash p_1(V(H))\times p_2(N_{{}_{x'_i}G_2}[V({}_{x'_i}H)])|$. And we divide subcase 2.3 into the following.

\noindent{\bf Subcase 2.3.1. }$|p_1(V(H))|\geq 5(k_1-1)+2$.

In this case, $|S|\geq \sum_{x'_{i}\in p_1(V(H))}|D_i|+|\cup^{l}_{i=1}(N_{G_1}(x'_i)\backslash p_1(V(H))\times p_2(N_{{}_{x'_i}G_2}[V({}_{x'_i}H)])|>\sum_{x'_{i}\in p_1(V(H))}|D_i|\geq |p_1(V(H))|\kappa(G_2)\geq (5(k_1-1)+2)k_2$. Therefore, $|S|>4k_1k_2-4$ when $k_1\leq k_2\leq 2k_1-1$ or $k_2\leq k_1\leq 2k_2-1$, $|S|>4k_1k_2+4k_1-2k_2-6$ when $k_2\geq 2k_1-1$, and $|S|>4k_1k_2-2k_1+4k_2-6$ when $k_1\geq 2k_2-1$.

\noindent{\bf Subcase 2.3.2. }$4\leq |p_1(V(H))|\leq 5(k_1-1)+1$.

Clearly, $H$ is connected, therefore $G_1[p_1(V(H))]$ is connected.

If $5\leq |p_1(V(H))|\leq 5(k_1-1)+1$, assume, without loss of generality, $G_1[\{x'_{1},x'_{2}, x'_{3},x'_{4},x'_{5}\}]$ is connected. Let $(x'_i,y'_i)\in V(H)$ for $i=1,2,3,4,5$. Since $girth(G_1)\geq 7$, $|N_{G_1}(p_1(V(H)))|\geq |N_{G_1}[\{x'_{1},x'_{2},x'_{3},x'_{4},x'_{5}\}]|-|p_1(V(H))|\geq 5k_1-3-l$. It implies that there are $5k_1-3-l$ neighbors of $V(H)$ in $\cup^{5}_{i=1}(N_{G_1}[x'_i]\backslash p_1(V(H)))\times \{y'_i\}$. These neighbors are denoted by $\Rmnum{1}$-type. Observe that $\cup^{l}_{i=1}(N_{G_1}[x'_i]\backslash p_1(V(H)))\times p_2(D_i)\subseteqq N_G(V(H))$. And $\sum^{5}_{i=1}|(N_{G_1}[x'_i]\backslash p_1(V(H)))\times p_2(D_i)|\geq (5k_1-3-l)k_2$. Therefore we have $(5k_1-3-l)k_2$ neighbors of $V(H)$ contained in $\cup^{5}_{i=1}(N_{G_1}[x'_i]\backslash p_1(V(H)))\times p_2(D_i)$. These neighbors are denoted by $\Rmnum{2}$-type. Since $girth(G_1)\geq 7$, $\Rmnum{1}$-type and $\Rmnum{2} $-type do not have the same element. Therefore there are at least $5k_1-3-l+(5k_1-3-l)k_2$ neighbors of $V(H)$ in $\cup^{l}_{i=1}(N_{G_1}(x'_i)\backslash p_1(V(H))\times p_2(N_{{}_{x'_i}G_2}[V({}_{x'_i}H)])$. Then $|S|\geq |N_{G}(V(H))|\geq \sum_{x'_{i}\in p_1(V(H))}|D_i|+|\cup^{l}_{i=1}(N_{G_1}(x'_i)\backslash p_1(V(H))\times p_2(N_{{}_{x'_i}G_2}[V({}_{x'_i}H)])|\geq lk_2+5k_1-3-l+(5k_1-3-l)k_2\geq (5k_1-3)k_2+1$. Therefore, $|S|>4k_1k_2-4$ when $k_1\leq k_2\leq 2k_1-1$ or $k_2\leq k_1\leq 2k_2-1$, $|S|>4k_1k_2+4k_1-2k_2-6$ when $k_2\geq 2k_1-1$, and $|S|>4k_1k_2-2k_1+4k_2-6$ when $k_1\geq 2k_2-1$.

If $|p_1(V(H))|=4$, assume, without loss of generality, $G_1[\{x'_{1}x'_{2}x'_{3}x'_{4}\}]$ is connected. Since $girth(G_1)\geq 7$, $|N_{G_1}(p_1(V(H)))|= |N_{G_1}(\{x'_{1},x'_{2},x'_{3},x'_{4}\})|\geq 4k_1-6$. This implies that there are $4k_1-6$ vertices of $\Rmnum{1}$-type and $(4k_1-6)k_2$ vertices of $\Rmnum{2}$-type. Therefore, there are at least $4k_1-6+(4k_1-6)k_2$ neighbors of $V(H)$ in $\cup^{l}_{i=1}(N_{G_1}(x'_i)\backslash p_1(V(H))\times p_2(N_{{}_{x'_i}G_2}[V({}_{x'_i}H)])$. Then $|S|\geq |N_{G}(V(H))|\geq \sum_{x'_{i}\in p_1(V(H))}|D_i|+|\cup^{l}_{i=1}(N_{G_1}(x'_i)\backslash p_1(V(H))\times p_2(N_{{}_{x'_i}G_2}[V({}_{x'_i}H)])|\geq 4k_2+ 4k_1-6+(4k_1-6)k_2= 4k_1k_2+4k_1-2k_2-6$. Therefore, $|S|\geq 4k_1k_2-4$ when $k_1\leq k_2\leq 2k_1-1$ or $k_2\leq k_1\leq 2k_2-1$, $|S|\geq 4k_1k_2+4k_1-2k_2-6$ when $k_2\geq 2k_1-1$, and $|S|>4k_1k_2-2k_1+4k_2-6$ when $k_1\geq 2k_2-1$.

According to the argument of subcase 2.3.1. and subcase 2.3.2., we can get that if any component of ${}_{x'_{i}}H$ contains at most 3 vertices for $i\in \{1,\cdots,l\}$ and $|p_1(V(H))|\geq 4$, then $|S|\geq 4k_1k_2-4$ when $k_1\leq k_2\leq 2k_1-1$ or $k_2\leq k_1\leq 2k_2-1$, $|S|\geq 4k_1k_2+4k_1-2k_2-6$ when $k_2\geq 2k_1-1$, and $|S|>4k_1k_2-2k_1+4k_2-6$ when $k_1\geq 2k_2-1$. Therefore, it remains to show that if any component of ${}_{x'_{i}}H$ contains at most 3 vertices for $i\in \{1,\cdots,l\}$ and $|p_1(V(H))|\leq 3 $, then $|S|\geq 4k_1k_2-4$ when $k_1\leq k_2\leq 2k_1-1$ or $k_2\leq k_1\leq 2k_2-1$, $|S|\geq 4k_1k_2+4k_1-2k_2-6$ when $k_2\geq 2k_1-1$, and $|S|\geq 4k_1k_2-2k_1+4k_2-6$ when $k_1\geq 2k_2-1$. Denote $H_y=G_{1y}\cap H$ for $y\in V(G_2)$, $p_2(V(H))=\{y'_1,\cdots,y'_h\}$ and $T_j=N_{G_{1y'_j}}(V(H_{y'_j}))$ for $j\in\{1,\cdots,h\}$. Since $|p_1(V(H)))|\leq 3$, $|V(H_{y'_j})|\leq 3$ for $j\in\{1,\cdots,h\}$. Doing similar arguments as subcase 2.3.1. and subcase 2.3.2., we can get that if $|p_2(V(H))|\geq 4$, then $|S|\geq 4k_1k_2-4$ when $k_1\leq k_2\leq 2k_1-1$ or $k_2\leq k_1\leq 2k_2-1$, $|S|> 4k_1k_2+4k_1-2k_2-6$ when $k_2\geq 2k_1-1$, and $|S|\geq 4k_1k_2-2k_1+4k_2-6$ when $k_1\geq 2k_2-1$. Since $S$ is a $\kappa_3$-$cut$ in $G$, it is impossible that $|p_1(V(H))|=1$ or $|p_2(V(H))|= 1$. This implies that we only need to show in the remaining proof that if $|p_1(V(H))|=2\text{ or } 3$ and $|p_2(V(H))|=2\text{ or } 3$ then $|S|\geq 4k_1k_2-4$ when $k_1\leq k_2\leq 2k_1-1$ or $k_2\leq k_1\leq 2k_2-1$, $|F_1|+|F_2|-|F_1\cap F_2|+|N_G(V(H))\backslash (F_1\cup F_2)|\geq 4k_1k_2+4k_1-2k_2-6$ when $k_2\geq 2k_1-1$, and $|F_1|+|F_2|-|F_1\cap F_2|+|N_G(V(H))\backslash (F_1\cup F_2)|\geq 4k_1k_2-2k_1+4k_2-6$ when $k_1\geq 2k_2-1$.

 Denote $F_1=\cup_{i=1}^{l}D_i$ and $F_2=\cup_{j=1}^{h}T_j$. Clearly, $|N_G(V(H))|=|F_1|+|F_2|-|F_1\cap F_2|+|N_G(V(H))\backslash (F_1\cup F_2)|$. Therefore, we only need to show that $|F_1|+|F_2|-|F_1\cap F_2|+|N_G(V(H))\backslash (F_1\cup F_2)|\geq 4k_1k_2-4$ when $k_1\leq k_2\leq 2k_1-1$ or $k_2\leq k_1\leq 2k_2-1$, $|S|\geq 4k_1k_2+4k_1-2k_2-6$ when $k_2\geq 2k_1-1$, and $|S|\geq 4k_1k_2-2k_1+4k_2-6$ when $k_1\geq 2k_2-1$. 

If $F_1\cap F_2=\emptyset$, then $p_2(V({}_{x'_{1}}H))=\cdots=p_2(V({}_{x'_{k}}H))$ and $p_1(V(H_{y'_{1}}))\cdots=p_1(V(H_{y'_{h}}))$. Since $\cup^{l}_{i=1} p_2(V({}_{x'_{i}}H))=p_2(V(H))$ and $\cup^{h}_{j=1} p_1(V(H_{y'_{j}}))=p_1(V(H))$, we can get $p_2(V({}_{x'_{i}}H)))=p_2(V(H))$ for $i\in\{1,\cdots,l\}$ and $p_1(V(H_{y'_{j}}))=p_1(V(H))$ for $j\in\{1,\cdots,h\}$. Therefore, $|F_1|+|F_2|=\sum^{l}_{i=1}|D_i|+\sum^{h}_{j=1}|T_j|\geq 2(2k_2-2)+2(2k_1-2)=4k_1+4k_2-8$. And $|N_G(V(H))\backslash (F_1\cup F_2)|=\sum_{i=1}^{l}|N_{G_1}[x'_{i}]\backslash p_1(V(H))\times D_i|\geq (2k_1-2)(2k_2-2)$. Thus $|N_G(V(H))|\geq |F_1|+|F_2|-|F_1\cap F_2|+|N_G(V(H))\backslash (F_1\cup F_2)|\geq 4k_1+4k_2-8+(2k_1-2)(2k_2-2)=4k_1k_2-4$. Therefore, $|S|\geq4k_1k_2-4$ when $k_1\leq k_2\leq 2k_1-1$ or $k_2\leq k_1\leq 2k_2-1$, $|S|\geq 4k_1k_2+4k_1-2k_2-6$ when $k_2\geq 2k_1-1$, and $|S|\geq 4k_1k_2-2k_1+4k_2-6$ when $k_1\geq 2k_2-1$.

If $F_1\cap F_2\neq\emptyset$, we construct a vertex set $R$ to count some elements of $N_G(V(H))\backslash (F_1\cup F_2)$ which are adjacent to $F_1\cap F_2$. Since $G_1[p_1(V(H)))]$ is connected, there is at least an edge $x'_{i_1}x'_{i_2}\in E(G_1[p_1(V(H)))])$ such that $p_2(V({}_{x_{i_1}}H))\neq p_2(V({}_{x_{i_2}}H))$. Let $(x'_{i_1},y'_{j_{i_1}})\in F_1\cap F_2$. If there is a vertex $y'\in N_{G_2}(y'_{j_{i_1}})\backslash p_2(V(H))$, then $(x'_{i_1},y')\in N_{{}_{x'_{i_1}}G_2}(V({}_{x_{i_2}}H))$(see Figure 2 for an illustration). Analogously, let $(x'_{i_2},y'_{j_{i_2}})\in F_1\cap F_2$. If there is a vertex $y''\in N_{G_2}(y'_{j_{i_2}})\backslash p_2(V(H))$, then $(x'_{i_2},y'')\in N_{{}_{x'_{i_2}}G_2}(V({}_{x_{i_1}}H))$. $R$ is defined as follows. 

\begin{align*}
R=\bigcup_{u\in F_1\cap F_2} p_1(u)\times (N_{G_2}(p_2(u))\backslash p_2(V(H))).
\end{align*}

\begin{figure}[htbp]
\centering
\includegraphics[width=8.5cm]{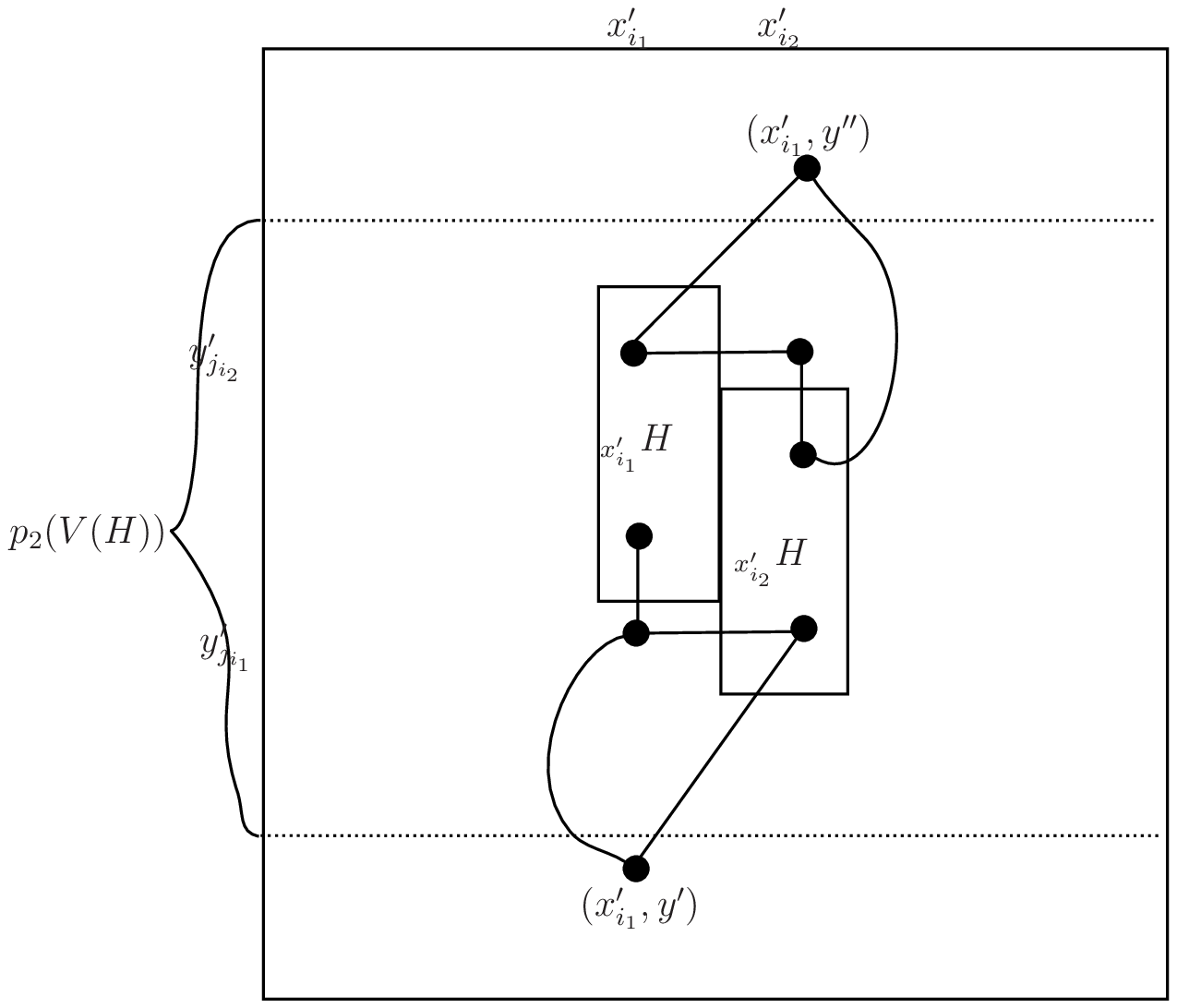}\\
\caption{An illustration of $(x'_{i_1},y')$ and $(x'_{i_2},y'')$.}
\end{figure}

And then we doing the argument according to $|p_1(V(H)|$ and $|p_2(V(H))|$ under $F_1\cap F_2\neq\emptyset$. 

$|p_1(V(H))|=|p_2(V(H))|=2$ was shown in the argument of $F_1\cap F_2=\emptyset$. If $|p_1(V(H))|=2$ and $|p_2(V(H))|= 3$, since $|V(H)|\geq 4$, we can get that $|F_1|=\sum_{i=1}^{2}|D_i|\geq 4k_2-4$, $|F_2|=\sum_{i=1}^{3}|T_j|\geq 4k_1-2$. Since $F_1\cap F_2\neq\emptyset$, $1\geq|F_1\cap F_2|\leq |p_1(V(H)\times p_2(V(H))|-|V(H)|\leq 2$. Observe that $|R|= \sum_{u\in F_1\cap F_2}|N_{G_2}(p_2(u))\backslash p_2(V(H))|\geq k_2-2$. Then $|N_G(V(H))\backslash (F_1\cup F_2)|\geq \sum_{i=1}^{2}|N_{G_1}[x'_{i}]\backslash p_1(V(H))\times D_i| +|R|\geq (k_1-1)(4k_2-4)+k_2-2$. Therefore, $|N_G(V(H))|\geq |F_1|+|F_2|-|F_1\cap F_2|+|N_G(V(H))\backslash (F_1\cup F_2)|\geq 4k_1+4k_2-8+(k_1-1)(4k_2-4)+k_2-2=4k_1k_2+k_2-6$. It implies that $|S|\geq4k_1k_2-4$ when $k_1\leq k_2\leq 2k_1-1$ or $k_2\leq k_1\leq 2k_2-1$, $|S|> 4k_1k_2+4k_1-2k_2-6$ when $k_2\geq 2k_1-1$, and $|S|\geq 4k_1k_2-2k_1+4k_2-6$ when $k_1\geq 2k_2-1$.

If $|p_1(V(H))|=3$ and $|p_2(V(H))|=2$. Since $|V(H)|\geq 4$, we can get that $|F_1|=\sum_{i=1}^{3}|D_i|\geq 4k_2-2$ and $|F_2|=\sum_{i=1}^{2}|T_j|\geq 4k_1-4$. Since $F_1\cap F_2\neq\emptyset$, $1\leq|F_1\cap F_2|\leq |p_1(V(H)\times p_2(V(H))|-|V(H)|\leq 2$. Observe that $|R|=\sum_{u\in F_1\cap F_2}|N_{G_2}(p_2(u))\backslash p_2(V(H))|\geq k_2-1$. Then $|N_G(V(H))\backslash (F_1\cup F_2)|= \sum_{i=1}^{3}|N_{G_1}[x'_{i}]\backslash p_1(V(H))\times D_i| +|R|=(k_1-2)\sum_{i=1}^{3}|D_i|+|D_1|+|D_3|+|R|\geq (k_1-2)(4k_2-2)+3k_2-1$. Therefore, $|N_G(V(H))|\geq |F_1|+|F_2|-|F_1\cap F_2|+|N_G(V(H))\backslash (F_1\cup F_2)|\geq 4k_1+4k_2-8+(k_1-2)(4k_2-2)+3k_2-1=4k_1k_2+2k_1-k_2-5$. It implies that $|S|\geq 4k_1k_2-4$ when $k_1\leq k_2\leq 2k_1-1$ or $k_2\leq k_1\leq 2k_2-1$, $|S|\geq 4k_1k_2+4k_1-2k_2-6$ when $k_2\geq 2k_1-1$, and $|S|> 4k_1k_2-2k_1+4k_2-6$ when $k_1\geq 2k_2-1$. 

If $|p_1(V(H))|=3$ and $|p_2(V(H))|=3$. When $|V(H)|\geq 5$,  we can get that $|F_1|=\sum_{i=1}^{3}|D_i|\geq 5k_2-4$ and $|F_2|=\sum_{i=1}^{3}|T_j|\geq 5k_1-4$. Since $F_1\cap F_2\neq\emptyset$, $1\leq|F_1\cap F_2|\leq |p_1(V(H)\times p_2(V(H))|-|V(H)|\leq 4$. Observe that $|R|=\sum_{u\in F_1\cap F_2}|N_{G_2}(p_2(u))\backslash p_2(V(H))|\geq k_2-2$. Then $|N_G(V(H))\backslash (F_1\cup F_2)|\geq \sum_{i=1}^{3}|N_{G_1}[x'_{i}]\backslash p_1(V(H))\times D_i| +|R|=(k_1-2)\sum_{i=1}^{3}|D_i|+|D_1|+|D_3|+|R|\geq (k_1-2)(5k_2-4)+3k_2-2$. Therefore, $|N_G(V(H))|\geq |F_1|+|F_2|-|F_1\cap F_2|+|N_G(V(H))\backslash (F_1\cup F_2)|\geq 5k_1+5k_2-12+(k_1-2)(5k_2-4)+2k_2=5k_1k_2+k_1-3k_2-4$. It implies that $|S|\geq4k_1k_2-4$ when $k_1\leq k_2\leq 2k_1-1$ or $k_2\leq k_1\leq 2k_2-1$, $|S|\geq 4k_1k_2+4k_1-2k_2-6$ when $k_2\geq 2k_1-1$, and $|S|\geq 4k_1k_2-2k_1+4k_2-6$ when $k_1\geq 2k_2-1$. When  $|V(H)|=4$, we can get that $|F_1|=\sum_{i=1}^{3}|D_i|\geq 4k_2-2$ and $|F_2|=\sum_{i=1}^{3}|T_j|\geq 4k_1-2$. Since $F_1\cap F_2\neq\emptyset$, $1\leq|F_1\cap F_2|\leq |p_1(V(H)\times p_2(V(H))|-|V(H)|\leq 5$. Then  $|R|=\sum_{u\in F_1\cap F_2}|N_{G_2}(p_2(u))\backslash p_2(V(H))|\geq k_2-1+k_2-2=2k_2-3$, and $|N_G(V(H))\backslash (F_1\cup F_2)|= \sum_{i=1}^{3}|N_{G_1}[x'_{i}]\backslash p_1(V(H))\times D_i| =(k_1-2)\sum_{i=1}^{3}|D_i|+|D_1|+|D_3|+|R|\geq (k_1-2)(4k_2-2)+4k_2-3$. Therefore,  $|N_G(V(H))|\geq|F_1|+|F_2|-|F_1\cap F_2|+|N_G(V(H))\backslash (F_1\cup F_2)|\geq 4k_1+4k_2-9+(k_1-2)(4k_2-2)+4k_2-3=4k_1k_2+2k_1-8$. It implies that $|S|\geq4k_1k_2-4$ when $k_1\leq k_2\leq 2k_1-1$ or $k_2\leq k_1\leq 2k_2-1$, $|S|\geq 4k_1k_2+4k_1-2k_2-6$ when $k_2\geq 2k_1-1$, and $|S|> 4k_1k_2-2k_1+4k_2-6$ when $k_1\geq 2k_2-1$.

Therefore, if $C_i$ is a component of ${}_{x'_{i}}H$ for $i\in \{1,\cdots,l\}$, then $|S|\geq$ min$\{4k_1k_2+4k_1-2k_2-6,4k_1k_2-2k_1+4k_2-6,4k_1k_2-4\}$ holds for all case of $|V(C_i)|$. Combining with case 1, Thus $\kappa_3(G)\geq$ min$\{4k_1k_2+4k_1-2k_2-6,4k_1k_2-2k_1+4k_2-6,4k_1k_2-4\}$, the proof is complete. $\Box$

\section{The Extra Conditional Fault-diagnosability of the Strong Product of Two Regular and Maximally Connected Graphs under PMC Model}

The process of identifying faulty vertices is called the diagnosis of the system. The diagnosability is the largest number of faulty vertices which can be identified by the system. The $PMC$ model using the matching theory in bipartite graphs was proposed by Preparata, Metze, and Chien \cite{Preparata} for dealing with the multiprocessor system's self-diagnosis. Under this model, assumed that each processor, regardless of the taster's faulty or non-faulty status, can test its neighboring processors through links between them, and the results of the test are either 0 or 1. When a fault-free tester evaluates a fault-free tested processor, then the result of the test is 0; when a fault-free tester evaluates a faulty tested processor, then the result of the test is 1. The collection of total outcomes which is done by testing all processors is called a syndrome. 

Let $R_1$ and $R_2$ be two vertex sets of $G$. The symmetric difference of $R_1$ and $R_2$ is $R_1 \Delta R_2$, and $E[R_1, R_2]$ express edge set in which every edge with one end in $R_1$ and the other end in $R_2$. If two different vertex sets $R_1$ and $R_2$ form different syndromes, then $R_1$ and $R_2$ are called distinguishable, otherwise, $R_1$ and $R_2$ are called indistinguishable.  Under $PMC$ model, the two different vertex sets $R_1$ and $R_2$ can be deemed a distinguishable-pair if and only if $E[V(G)-R_1 \cup R_2, R_1 \Delta R_2]\neq\emptyset$ \cite{Lai}.

The classical diagnosis strategy ignores the fact that it is extremely unlikely that all neighboring processors of a processor can fail at the same time. Therefore, Lai \cite{Lai} developed a new measure of diagnosability, named conditional diagnosability of a system under the $PMC$ model, which assumes not all adjacent processors of one processor may be faulty simultaneously. Then Zhang and Yang\cite{Zhang1} introduced the $g$-$extra$ conditional fault-diagnosability. A system designed by $G$ is $g$-$extra$ conditionally $t$-fault-diagnosable if and only if for each pair of distinct $g$-$extra$ $cut$ $R_1, R_2\subseteq V(G)$ such that $|R_1| \leq t,|R_2| \leq t$ and $R_1$ and $R_2$ are distinguishable. The $g$-$extra$ conditional fault-diagnosability of $G$, is the largest value of $t$ such that $G$ is $g$-$extra$ conditionally $t$-fault-diagnosable. Clearly, The $g$-$extra$ connectivity and the $g$-$extra$ conditional fault-diagnosability have an intimate association.
For more results on $g$-$extra$ conditional fault-diagnosability of different networks please refer to $\cite{Liu, Lin, Yuan}$.

We use $\tilde{t_g^p}(G)$ to denote the $g$-$extra$ conditional fault-diagnosability of $G$ in the $PMC$ model. Lin, Xu, Chen, and Hsieh \cite{Lin1} established a relationship between $\kappa_g(G)$ and $\tilde{t_g^p}(G)$ for a $k$-regular graph $G$ with $1 \leq g \leq k$, under $PMC$ model.

\begin{theorem}\cite{Lin1}\label{9}
Given a general $k$-regular graph $G$ and an integer $g$ with $0 \leq g \leq k$. If $|V(G)| \geq 2[\kappa_g(G)+g]+1$ and there exists a connected subgraph $\tilde{S_h}$ such that $|\tilde{S_h}|=g+1$ and $N_G(V(\tilde{S_h}))$ be the minimum $g$-$extra$ $cut$ of $G$, then $\tilde{t_g^p}(G)=\kappa_g(G)+g$.
\end{theorem}

For the sake of simplicity, denote
$M_1(g,k_1,k_2)=(g+1)(k_1k_2-k_1-k_2)+2(g+1)($min$\{k_1,k_2\}-1)+2($max$\{k_1,k_2\}-1)+4$, $M_2(g,k_1,k_2)=(g+1)(k_1k_2-k_1-k_2)+2\lceil\sqrt{g+1}\,\rceil($min$\{k_1,k_2\}-1)+2\lceil\frac{g+1}{\lceil\sqrt{g+1}\,\rceil}\rceil($max$\{k_1,k_2\}-1)+4$ and $M(g,k_1,k_2)=$min$\{M_1(g,k_1,k_2),M_2(g,k_1,k_2)\}$. We can combine $Theorem\,3.1$ to $3.3$ into the following result directly.

\begin{theorem}\label{11}
Let $g\ (\leq 3)$ be a non-negative integer, and $G_i$ be a $k_i\ (\geq 2)$-regular and maximally connected graph with $girth(G_i)\geq g+4$ for $i=1,2$. Then $\kappa_{g}(G_1\boxtimes G_2)= M(g,k_1,k_2)$.
\end{theorem}

According to argument of section 3, we can find a connected subgraph $H$ in $G_1\boxtimes G_2$ such that $|V(H)|=g+1$ and $N_{G_1\boxtimes G_2}(V(H))$ be a $\kappa_g$-cut of $G_1\boxtimes G_2$ for $g=1,2,3$, where $G_i$ be a $k_i\ (\geq 2)$-regular and maximally connected graph with $girth(G_i)\geq g+4$ for $i=1,2$. Thus, by $Theorem\,\ref{9}$, we obtain the following corollary directly.

\begin{corollary}
Let $g\ (\leq 3)$ be a non-negative integer, and $G_i$ be a $k_i\ (\geq 2)$-regular and maximally connected graph with $girth(G_i)\geq g+4$ for $i=1,2$. Then $\tilde{t_g^p}(G_1\boxtimes G_2)=M(g,k_1,k_2)+g$.
\end{corollary}

\section{Concluding Remarks}
Graph products are used to construct large graphs from small ones. Strong product is one of the most studied four graph products. As a generalization of traditional connectivity, $g$-$extra$connectivity can be seen as a refined parameter to measure the reliability of interconnection networks. There is no polynomial-time algorithm to compute the $g\ (\geq1)$-extra connectivity for a general graph. In this paper, we determined the $g\ (\leq 3)$-extra connectivity of the strong product of two maximally connected regular graphs whose girth is no less than $g+4$. Furthermore, under the $PMC$ model, we got the $g\ (\leq 3)$-extra conditional fault-diagnosability of them. In future work, we would like to investigate the $g$-$extra$ connectivity and $g$-$extra$ conditional fault-diagnosability of more general graphs.



\end{sloppypar}
\end{document}